\documentclass[10pt,draft,a4paper]{amsart}

\usepackage{amsmath}
\usepackage{amssymb}
\usepackage{amsfonts}
\usepackage{graphicx}
\numberwithin{equation}{section}
\newtheorem{theorem}{Theorem}[section]

\newtheorem{lemma}[theorem]{Lemma}

\theoremstyle{remark}

\theoremstyle{definition}

\theoremstyle{remark}

\theoremstyle{remark}
\newtheorem{remark}{Remark}[section]

\begin{document}

\newcommand{\spt}{\,\mathrm{supp}\,}
\newcommand{\xxi}{\langle\xi\rangle}
\newcommand{\xx}{\langle x\rangle}
\newcommand{\yy}{\langle y\rangle}
\newcommand{\mm}{\langle \mu\rangle}
\newcommand{\dint}{\int\!\!\int}
\newcommand{\triple}[1]{{|\!|\!|#1|\!|\!|}}
\newcommand{\Sph}{{\mathbb S}}

\title[Schr\"odinger with point interaction]%
{Dispersive estimate for the Schr\"odinger equation
with point interactions}

\begin{abstract}
We consider the Schr\"{o}dinger operator in $\mathbb{R}^{3}$ with $N$
point interactions
placed at $Y=(y_{1},\ldots ,y_{N})$, $y_{j} \in \mathbb{R}^{3}$, of
strength $\alpha=(\alpha_{1}, \ldots ,\alpha_{N}) $, $\alpha_{j}\in
\mathbb{R}$. Exploiting the spectral theorem and the rather explicit
expression for the resolvent we prove the (weighted) dispersive estimate
for the corresponding Schr\"{o}dinger flow.

\noindent
In the special case $N=1$ the proof is directly obtained from the
unitary group which is known in closed form.

\end{abstract}

\date{\today} 

\author{Piero D'Ancona}
\address{Piero D'Ancona:
Unversit\`a di Roma ``La Sapienza'',
Dipartimento di Matematica,
Piazzale A.~Moro 2, I-00185 Roma, Italy}
\email{dancona@mat.uniroma1.it}

\author{Vittoria Pierfelice}
\address{Vittoria Pierfelice:
Universit\`a di Pisa,
Dipartimento di Matematica,
Via Buonarroti 2, I-56127 Pisa, Italy}
\email{pierfelice@dm.unipi.it}

\author{Alessandro Teta}
\address{Alessandro Teta:
Universit\`a di L'Aquila,
Dipartimento di Matematica Pura e Applicata,
Via Vetoio, Loc. Coppito, I-67100 L'Aquila, Italy}
\email{teta@univaq.it}

\maketitle

\section{Introduction}\label{sec.introd}

In this paper we study the time dependent Schr\"{o}dinger equation
in $\mathbb{R}^{3}$ with a finite number of point interactions and, in
particular, we shall prove a dispersive estimate for the solution.

\noindent
At a formal level the Schr\"{o}dinger operator with point interactions
can be written as

\begin{equation}\label{formal}
H= -\Delta+ \sum_{j=1}^{N} \mu_{j} \delta_{y_{j}}
\end{equation}

\noindent
where $\delta_{y_{j}}$ is the Dirac measure placed at $y_{j} \in
\mathbb{R}^{3}$ and the parameters $\mu_{j}$ are coupling constants.
As a matter of fact the Dirac measure in $\mathbb{R}^{3}$ is not
a small perturbation of the Laplacian, even in the sense of quadratic
forms.

\noindent
As a consequence a self-adjoint Hamiltonian in
$L^{2}(\mathbb{R}^{3})$ cannot be defined as the sum
of a kinetic plus an interaction part and this explains why \eqref{formal} is
only a formal expression.

\noindent
In order to obtain a rigorous counterpart of \eqref{formal} one
considers the following restriction of the free Laplacian

\begin{equation}\label{restr}
\hat{H} = - \Delta , \;\;\;\;\; D(\hat{H})
=C_{0}^{\infty}(\mathbb{R}^{3} \setminus Y)
\end{equation}

\noindent
where $Y=(y_{1}, \ldots ,y_{N})$. The operator \eqref{restr} is
symmetric but not self-adjoint in $L^{2}(\mathbb{R}^{3})$ and,
obviously, one possible self-adjoint extension is trivial, i.e. it
coincides with the free Laplacian
$H_{0}=- \Delta$, $D(H_{0})=H^{2}(\mathbb{R}^{3})$.

\noindent
Using the theory of self-adjoint extensions of symmetric operators,
developed by von Neumann and Krein, one can show that the operator
\eqref{restr} has $N^{2}$ (non trivial) self-adjoint extensions which,
by definition, are all the possible Schr\"{o}dinger operators with
point interactions at $Y$ (for a comprehensive treatment we refer to
the monograph \cite{Albeverio88-solvablemod}).

\noindent
Any such extension can be considered as a Laplace operator with a singular boundary
condition satisfied at each point $y_{j} \in Y$.

\noindent
In the following we
shall only consider the case of local boundary conditions which are
more relevant from the physical point of view. More precisely, we
shall restrict to the self-adjoint extensions $H_{\alpha,Y}$ parametrized by
$\alpha =(\alpha_{1}, \ldots,\alpha_{N})$, $\alpha_{j} \in
\mathbb{R}$, and corresponding to the singular boundary condition at
$Y$

\begin{equation}\label{bc}
\lim_{r_j \rightarrow 0} \left[ \frac{\partial(r_j u)}{\partial r_j} - 4 \pi
\alpha_j (r_j u) \right]=0, \;\;\;\;r_j=|x-y_j|,\;\;\;\;j=1,...,N
\end{equation}

\noindent
In the next section we shall give the precise definition and the main
properties of the Hamiltonian $H_{\alpha,Y}$.

\noindent
We only notice that the domain of $H_{\alpha,Y}$ contains functions
with a singularity of the type $|x-y_{j}|^{-1}$ at each point $y_{j}
\in Y$ and this explains why in the boundary condition \eqref{bc} the
behaviour of the function $u$ near $y_{j}$ must be regularized.

\noindent
We also recall that the physical meaning of the parameters
$\alpha_{j}$ is connected with the scattering length of the scatterer
in $y_{j}$ and their relation with the formal coupling constants
$\mu_{j}$ introduced in \eqref{formal} can be understood via a
suitable renormalization procedure (see e.g.
\cite{Albeverio88-solvablemod}).

\noindent
In the following we shall be concerned with the unitary evolution
group generated by $H_{\alpha,Y}$

\begin{equation}\label{eq.propa}
e^{itH_{\alpha,Y}}f
\end{equation}

\noindent
which gives the solution of the Cauchy problem

\begin{equation}
iu_{t}+ H_{\alpha,Y}u=0,
\qquad u(0,x)=f(x).
\end{equation}

\noindent
It is well known that for the free Schr\"{o}dinger group in
dimension three the following
dispersive estimate holds

\begin{equation}\label{eq.dispclass}
\left\|e^{-it\Delta}f\right\|_{L^{\infty}(\mathbb{R}^{3})}
\leq \frac C{t^{3/2}}\|f\|_{L^{1}(\mathbb{R}^{3})}
\end{equation}

\noindent
The estimate \eqref{eq.dispclass} can be generalized to
smooth perturbation of the Laplacian, i.e. the l.h.s. can be replaced
by $P_{ac} e^{it(-\Delta +V)}f$, where $V$ is a smooth
real-valued potential satisfying some decay condition at infinity,
$P_{ac}$ is the projection onto the absolutely continuous spectrum
of $-\Delta +V$ and, moreover, it is assumed that zero is not an eigenvalue nor a
resonance for $-\Delta +V$ (see e.g.
\cite{GoldbergSchlag}, \cite{RodnianskiSchlag04-disp}
\cite{Yajima95-waveopN}, \cite{Yajima95-waveopNeven}
and the references therein).

\noindent
Here we want to extend the validity of the dispersive estimate to the
case of the Schr\"{o}dinger equation with point interactions.

\noindent
Due to the presence of the (unavoidable) singularity at the points where the interaction is placed,
we are forced to introduce the following weight function

\begin{equation}\label{eq.weight}
w(x)=\sum_{j=1}^{N}
\left(
1+\frac1{|x-y_j|}
\right).
\end{equation}

\noindent
Moreover, for $z \in \mathbb{C}$, we define the matrix

\begin{equation}\label{eq.Gamma1}
[\Gamma_{\alpha,Y}(z)]_{j,\ell}= \left[
\left( \alpha_{j}-\frac{iz}{4\pi} \right) \delta_{j,\ell} - \tilde{G}_z(y_j-y_\ell)
\right]_{j,\ell=1}^N,
\end{equation}
with
\begin{equation}\label{eq.Gamma2}
\tilde{G}_z(x)=
\begin{cases}\displaystyle
\frac{e^{iz|x|}}{4\pi|x|} \quad \;\; \;x\neq 0,\\
0 \qquad \qquad x=0.
\end{cases}
\end{equation}


\noindent
The role of the matrix $\Gamma_{\alpha,Y}(z)$ for the characterization
of the main properties of $H_{\alpha,Y}$ will be explained in the next
section. Here we only notice that our basic assumption is the
invertibility of $\Gamma_{\alpha,Y}(\mu)$, for $\mu \in [0,+
\infty)$. Such assumption in particular implies that zero is not an
eigenvalue nor a resonance for $H_{\alpha,Y}$.

\noindent
In the special case $N=1$, i.e. for the Schr\"{o}dinger operator
$H_{\alpha,y}$ with a
point interaction in $y \in \mathbb{R}^{3}$ of strength $\alpha \in
\mathbb{R}$, the assumption simply means $\alpha \neq
0$, which is precisely the condition for the absence of a zero-energy
resonance.

\noindent
Our main result is the following theorem.

\begin{theorem}\label{th.1}
Assume that the matrix $\Gamma_{\alpha,Y}(\mu)$ is invertible
for $\mu\in[0,+\infty)$ with a locally bounded inverse. Then
the following dispersive estimate holds
\begin{equation}\label{eq.dispgrossa}
\left\|w^{-1}e^{itH_{\alpha,Y}}P_{ac}f\right\|
_{L^{\infty}(\mathbb{R}^{3})}
\leq
\frac C{t^{3/2}}\|w\cdot f\|_{L^{1}(\mathbb{R}^{3})}
\end{equation}
for any $f \in L^{2}(\mathbb{R}^{3})$ such that $w \cdot f \in
L^{1}(\mathbb{R}^{3})$.

\noindent
In the special case $N=1$, estimate \eqref{eq.dispgrossa}
holds for all $\alpha\neq0$; moreover, when $\alpha>0$
the projection
$P_{ac}$ can be replaced by the identity.
Finally, in the resonant case $\alpha=0$ we have
the slower decay estimate
\begin{equation}\label{eq.disp2}
\left\|w^{-1}e^{itH_{0,y}}f\right\|_{L^{\infty}(\mathbb{R}^{3})}
\leq
\frac C{t^{1/2}}\|w\cdot f\|_{L^{1}(\mathbb{R}^{3})}
\end{equation}
\end{theorem}

\noindent
We remark that the result for $N=1$ is more detailed due
to the fact that the unitary propagator is explicitely known (see
e.g. \cite{ScarlattiTeta90}, \cite{AlbeverioBrzezniakDabrowski}) and then
the estimate can be obtained by a straightforward computation.

\noindent
The rest of the paper is organized as follows.

\noindent
In section 2 we recall the precise definition and the main properties
of the Hamiltonian $H_{\alpha,Y}$.

\noindent
In section 3 we give the proof of the theorem in the special case
$N=1$, exploiting the knowledge of the kernel of the Schr\"{o}dinger
group.

\noindent
In section 4 we give the proof in the general case using the
spectral calculus and the explicit expression for the resolvent of $H_{\alpha,Y}$.

\section{The Schr\"odinger operator with point interactions}\label{sec.delta}

In this section we review the definition and some basic properties of the Schr\"odinger operator with point interactions,
referring to \cite{Albeverio88-solvablemod} for more details.
We start considering the perturbation of $-\Delta$ by
a zero-range (singular) potential
supported at the origin of $\mathbb{R}^3$. As we pointed out in the
introduction, such perturbation is by definiton a non trivial self-adjoint
extension of the symmetric operator \eqref{restr} when $N=1$.

We recall that the trivial extension $-\Delta$ is the (self-adjoint) Laplace operator in
$L^2(\mathbb{R}^3)$
(with Lebesgue measure),
with domain
$D(-\Delta)= H^2(\mathbb{R}^3)= W^{2,2}(\mathbb{R}^3)$,
where $W^{p,k}(\mathbb{R}^{3})$ denotes the Sobolev space
of all functions belonging to
$L^p(\mathbb{R}^3)$
whose weak derivatives
of order smaller or equal to $k$ belong to $L^p(\mathbb{R}^3)$.

Using the theory of self-adjoint extensions, one can show that the
non trivial self-adjoint extensions $H_{\alpha}$, parametrized by
$\alpha \in \mathbb{R}$, can be completely
characterized.

Indeed, the domain $\mathcal{D}(H_\alpha)$ consists of all elements
$\psi \in L^{2}(\mathbb{R}^{3})$ of the type
\begin{equation}\label{deco}
\psi(x) = \phi_z(x) + \left(\alpha-\frac{i z}{4\pi}\right)^{-1} \phi_z(0)
\frac{e^{- i z|x|}}{4\pi|x|}
\end{equation}
where $\phi_z \in H^2(\mathbb{R}^3)$ and ${\Im}z>0$.
The decomposition \eqref{deco} is unique and with $\psi \in
\mathcal{D}(H_\alpha)$ one has
\begin{equation}
(H_\alpha - z^2) \psi= (-\Delta - z^2) \phi_z.
\end{equation}

We notice that the trivial extension is recovered in the limit
$\alpha \rightarrow \infty$.

A remarkable property of the Hamiltonian $H_\alpha$ is that the integral kernel
of the resolvent $R_{\alpha}(z)=(H_{\alpha}- z)^{-1}, \; {\Im}z > 0,$ can be
explicitely computed by the krein's formula. In fact

\begin{equation}\label{eq.kernR1}
(R_{\alpha}(z^2)f)(x)=
(R_{0}( z^2)f)(x) +
\left(\alpha-\frac{i z}{4\pi}\right)^{-1}
\frac{e^{- i z|x|}}{4\pi|x|}
\int_{\mathbb{R}^{3}}
\frac{e^{ i z|y|}}{4\pi|y|}
f(y)dy.
\end{equation}

where the free resolvent is given by
\begin{equation*}
(R_{0}(z^{2})f)(x)=
\int_{\mathbb{R}^{3}}
\frac{e^{i z|x-y|}}{4\pi|x-y|}
f(y)dy,
\end{equation*}
By an elementary computation we obtain that for any $\lambda\in\mathbb{R}$ and $\varepsilon>0$
\begin{equation}\label{freepm}
(R_0(\lambda \pm i\varepsilon)f)(x)=
\int_{\mathbb{R}^{3}} \frac{e^{\pm i\sqrt{\lambda_\varepsilon} |x-y|}}{4 \pi |x-y|}
e^{-\varepsilon|x-y|/2\sqrt{\lambda_{\varepsilon}}} f(y) dy,
\end{equation}
where
\begin{equation}\label{laep}
\lambda_\varepsilon =\frac{\lambda+(\lambda^{2}+\varepsilon^{2})^{1/2}}{2}>0.
\end{equation}

From the last inequality it is easy to derive the limiting absorption
principle for the free resolvent, i.e. for
$\lambda>0$
\begin{equation}\label{eq.kernR0}
(R_{0}(\lambda+ i0)f)(x)=
\int_{\mathbb{R}^{3}}
\frac{e^{i \sqrt \lambda|x-y|}}{4\pi|x-y|}
f(y)dy.
\end{equation}

We can also deduce the limiting absorption principle for the
resolvent \eqref{eq.kernR1}
\begin{equation}\label{eq.kernR}
(R_{\alpha}(\lambda+ i0)f)(x)=
(R_{0}(\lambda+ i0)f)(x)+
\left(\alpha-\frac{i\sqrt\lambda}{4\pi}\right)^{-1}
\frac{e^{- i \sqrt\lambda|x|}}{4\pi|x|}
\int_{\mathbb{R}^{3}}
\frac{e^{ i \sqrt\lambda|y|}}{4\pi|y|}
f(y)dy.
\end{equation}

The spectral properties of $H_\alpha$ are easily derived from
\eqref{eq.kernR1}.
The continuous spectrum of $H_\alpha$ is purely absolutely continuous and covers the nonnegative real axis
$[0,\infty)$
while the point spectrum is empty if $\alpha\geq 0$ and $\{(-4\pi \alpha)^2\}$ if $\alpha<0$.
The normalized eigenfunction associated with the only negative eigenvalue is given by
\begin{equation}\label{eq.eigf}
\Psi_{\alpha}(x)=\sqrt{-2\alpha}\
\frac{e^{4\pi\alpha|x|}}{|x|},\qquad
\alpha<0.
\end{equation}
For $\alpha=0$ the Hamiltonian has a zero-energy resonance.

In analogous way we introduce the properties of the Schr\"odinger operator $H_{\alpha,Y}$ with point
interactions located at
$Y=(y_1,\ldots,y_N)$ with strength $\alpha=(\alpha_1,\ldots,\alpha_N)$.

The domain $\mathcal{D}(H_{\alpha,Y})$ is the set of all functions
$\psi \in L^{2}(\mathbb{R}^{3})$ of following type

\begin{equation}\label{deco2}
\psi(x) = \phi_z(x) + \sum_{j=1}^N a_j
\frac{e^{- i z|x - y_j|}}{4\pi|x - y_j|}, \;\;\;\;\; x \in
\mathbb{R}^3 \setminus Y,
\end{equation}
where $$a_j= \sum_{l=1}^{N}
[\Gamma_{\alpha,Y}(z)]^{-1}_{j,l} \phi_z(y_l), \quad j=1, \ldots,N,$$
$\phi_z \in H^2(\mathbb{R}^3)$,\; ${\Im}z>0$ and the matrix
$\Gamma_{\alpha,Y}(z)$ has been defined in \eqref{eq.Gamma1}.
The decomposition \eqref{deco2} is unique and with $\psi \in \mathcal{D}(H_{\alpha,Y})$ we obtain
\begin{equation}
(H_{\alpha,Y} - z^2) \psi= (-\Delta - z^2) \phi_z.
\end{equation}
It is an easy computation to verify that each element of $\mathcal{D}(H_{\alpha,Y})$
satisfies the local boundary condition \eqref{bc}.

Also in the $N$-centers case the resolvent of these operators is explicitly given by Krein's formula.
In fact $R_{\alpha,Y}(z)=(H_{\alpha,Y}- z)^{-1}$ is given by the following formula
for $\Im z>0$
\begin{equation}\label{eq.kernR2}
(R_{\alpha,Y}(z^{2})f)(x) = (R_{0}(z^{2})f)(x) +
\sum_{j,\ell=1}^{N}
[\Gamma_{\alpha,Y}(z)]^{-1}_{j,\ell}
\frac{e^{-i z|x-y_{j}|}}{4\pi|x-y_{j}|}
\int_{\mathbb{R}^{3}}
\frac{e^{i z|y-y_{\ell}|}}{4\pi|y-y_{\ell}|}
f(y)dy
\end{equation}
Moreover we deduce the validity of the limiting absorption principle,
i.e. for
$ \lambda>0$
\begin{equation}\label{eq.kernRr}
(R_{\alpha,Y}(\lambda+i0)f)(x)=(R_{0}(\lambda+i0)f)(x)+
\sum_{j,\ell=1}^{N}
[\Gamma_{\alpha,Y}(i\sqrt \lambda)]^{-1}_{j,\ell}\
\frac{e^{-i\sqrt \lambda |x-y_{j}|}}{4\pi|x-y_{j}|}
\int_{\mathbb{R}^{3}}
\frac{e^{i\sqrt \lambda |y-y_{\ell}|}}{4\pi|y-y_{\ell}|}
f(y)dy
\end{equation}

We conclude summarizing the spectral properties of $H_{\alpha,Y}$.
The essential spectrum of $H_\alpha$ is purely absolutely continuous and coincides with the
real axis $[0,\infty)$, indeed the singular continuous spectrum is empty.
Moreover, the operator $H_{\alpha,Y}$ has the point spectrum included
in $ (0, -\infty)$, i.e. there are no positive
embedded eigenvalues; in particular $H_{\alpha,Y}$ has at most $N$ (negative) eigenvalues counting
multiplicity.
In addition, the eigenvalues can be determined
computing the zeros
of the determinant of an $N\times N$ matrix, i.e.
\begin{equation*}
z^{2}\in\sigma_p(H_{\alpha,Y})\iff
\det\Gamma_{\alpha,Y}(z)=0,
\end{equation*}
and the multiplicity of eigenvalue $z^2$ equals the multiplicity of the eigenvalue zero of the matrix
$\Gamma_{\alpha,Y}(z)$.

\section{Proof of the result for $N=1$}\label{sec.proofth.1a}


In the case of a single point interaction considered
in in the second part of Theorem \ref{th.1}, the proof is quite easy.
We believe it is of some interest, both because the
proof is much easier, and it is possible to
give a complete description of the dispersive behaviour
of the propagator, including the case when
a resonance at zero occurs.
The main advantage
is that, as proved in \cite{ScarlattiTeta90}
(see also \cite{AlbeverioBrzezniakDabrowski}),
it is possible to give an explicit representation
of the Schr\"odinger propagator. If we denote by
\begin{equation}\label{eq.prop}
S(x;t)=\frac{e^{|x|^{2}/4it}}{(4\pi it)^{3/2}},\qquad
t>0,\quad x,y\in\mathbb{R}^{3}
\end{equation}
the free Schr\"odinger propagator in $\mathbb{R}^{3}$, i.e.
\begin{equation*}
(e^{-it\Delta}f)(x)=S(x;t)*_{x}f(x),
\end{equation*}
the kernel of the propagator of $e^{itH_{\alpha}}$ can be written as follows:

1) for $\alpha>0$
\begin{equation}\label{eq.prop1}
S_{\alpha}(x,y;t)=
S(x-y;t)
+\frac 1{|x||y|}\int_{0}^{\infty}
e^{-4\pi\alpha s}(s+|x|+|y|)
S(s+|x|+|y|;t)ds
\end{equation}

2) for $\alpha=0$
\begin{equation}\label{eq.prop2}
S_{0}(x,y;t)=
S(x-y;t)
+\frac {2it}{|x||y|}
S(|x|+|y|;t)
\end{equation}

3) for $\alpha<0$

\begin{equation}\label{eq.prop3}
S_{\alpha}(x,y;t)=
S(x-y;t)
+\Psi_{\alpha}(x)\Psi_{\alpha}(y)e^{it(4\pi\alpha)^{2}}
+\frac 1{|x||y|}\int_{0}^{\infty}
e^{4\pi\alpha s}(s-|x|-|y|)
S(s-|x|-|y|;t)ds
\end{equation}
Then we have the representation
\begin{equation*}
(e^{itH_{\alpha}}f)(x)=
\int_{\mathbb{R}^{3}}S_{\alpha}(x,y;t)f(y)dy.
\end{equation*}
Notice that in the case $\alpha<0$
we can easily distinguish the continuous part from the
standing wave: indeed, the continuous part can be expressed as
\begin{equation}\label{eq.Pac}
(P_{ac} e^{itH_{\alpha}}f)(x)=
(e^{-it\Delta}f)(x)+
\int_{\mathbb{R}^{3}}
\int_{0}^{\infty}\frac 1{|x||y|}
e^{4\pi\alpha s}(s-|x|-|y|)
S(s-|x|-|y|;t)dsdy.
\end{equation}
The standing wave is the function
\begin{equation*}
C_{1}\Psi_{\alpha}(x)e^{it(4\pi\alpha)^{2}},\qquad
C_{1}=\int_{\mathbb{R}^{3}}
\Psi_{\alpha}(y)f(y)dy
\end{equation*}
which can be written explicitly using \eqref{eq.eigf} as
\begin{equation}\label{eq.stw}
e^{it(4\pi\alpha)^{2}}(-2\alpha)
\frac{e^{4\pi\alpha|x|}}{|x|}
\int_{\mathbb{R}^{3}}
\frac{e^{4\pi\alpha|y|}}{|y|}f(y)dy.
\end{equation}

Consider the case $\alpha>0$ first. From the trivial estimate
\begin{equation}\label{eq.triv1}
|S(s+|x|+|y|;t)|\leq\frac C{t^{3/2}}
\end{equation}
and
\begin{equation}\label{eq.triv2}
\int_{0}^{\infty}
e^{-4\pi\alpha s}(s+|x|+|y|)ds\leq
C(1+|x|+|y|)
\end{equation}
and using also the standard dispersive estimate for the free
solution, we obtain immediately the estimate
\begin{equation*}
\left|\int_{\mathbb{R}^{3}} S_{\alpha}(x,y;t)f(y)dy\right|\leq
\frac C{t^{3/2}}
\int_{\mathbb{R}^{3}} \left(1+\frac{1+|x|+|y|}{|x||y|}\right)|f(y)|dy.
\end{equation*}
Since we have
\begin{equation*}
1+\frac{1+|x|+|y|}{|x||y|}=
\left(1+\frac1{|x|}\right) \left(1+\frac1{|y|}\right)=
w(x)w(y)
\end{equation*}
we easily conclude that
\begin{equation*}
|(w^{-1}e^{itH_{\alpha}}f)(x)|\leq
\frac C{t^{3/2}}
\int_{\mathbb{R}^{3}} w(y)|f(y)|dy
\end{equation*}
as claimed.

In the case $\alpha<0$, we can estimate the continuous part
\eqref{eq.Pac} exactly in the same way, obtaining
\begin{equation*}
|(w^{-1}P_{ac}e^{itH_{\alpha}}f)(x)|\leq
\frac C{t^{3/2}}
\int_{\mathbb{R}^{3}} w(y)|f(y)|dy.
\end{equation*}

Finally, when $\alpha=0$ we can see the effect of the resonance at
zero in the slower rate of decay: from \eqref{eq.prop2},
proceeding exactly as before, we obtain
\begin{equation*}
|(w^{-1}e^{itH_{\alpha}}f)(x)|\leq
\frac C{t^{1/2}}
\int_{\mathbb{R}^{3}} w(y)|f(y)|dy
\end{equation*}
since an additional factor $t$ is present in the propagator.

\section{Proof of the result for $N>1$}\label{sec.proofth.1b}

We shall now prove the main part of Theorem \ref{th.1} concerning the
general case $N>1$, where we do not
have an explicit formula for the propagator. Thus we resort to
the spectral calculus and we represent the (continuous part of the)
solution as follows:
\begin{equation}\label{eq.spec}
P_{ac}e^{itH_{\alpha,Y}}f=\int_{0}^{+\infty}
e^{it\lambda}\Im R_{\alpha,Y}(\lambda+i0)fd\lambda.
\end{equation}
Recall now the integral expression \eqref{eq.kernRr}. The first term
reproduces the free solution $e^{-it\Delta}f$ for which we already know
the decay estimate.
Thus we are reduced to estimate the integrals
\begin{equation}\label{eq.I}
I=
\int_{0}^{+\infty}\int_{\mathbb{R}^{3}}
e^{it\lambda}
\Im
\left(
e^{-i\sqrt \lambda |x-y_{j}|}
e^{i\sqrt \lambda |y-y_{\ell}|}
c_{j\ell}(\sqrt\lambda)
\right)
\frac{f(y)}{|x-y_{j}||y-y_{\ell}|}
dyd\lambda
\end{equation}
where
\begin{equation}\label{eq.cjl}
c_{j\ell}(\mu)=
(4\pi)^{-2}
[\Gamma_{\alpha,Y}(\mu)]^{-1}_{j,\ell}
\qquad
\text{for $\mu\geq0$,\ \ $j,\ell=1,\dots,N$.}
\end{equation}
To make the following computation rigorous,
we introduce a cutoff function.
Let $\psi(s)$ be a nonnegative function in
$C^{\infty}_{0}(\mathbb{R})$ equal to 1 on $[0,1]$ and vanishing
on $[2,\infty)$; it is not restrictive to
assume that $\psi$ is an \emph{even} function
$\psi(-s)=\psi(s)$.
Then we can approximate $I$ with
\begin{equation}\label{eq.IM}
I_{M}=
\int_{0}^{+\infty}\int_{\mathbb{R}^{3}}
e^{it\lambda}
\Im
\left(
e^{-i\sqrt \lambda |x-y_{j}|}
e^{i\sqrt \lambda |y-y_{\ell}|}
c_{j\ell}(\sqrt\lambda)
\right)
\frac{f(y)}{|x-y_{j}||y-y_{\ell}|}
\psi(\sqrt\lambda/M)
dyd\lambda
\end{equation}
as $M\to+\infty$. After the change of variables $\lambda=\mu^{2}$
this can be written also as
\begin{equation}\label{eq.IMbis}
I_{M}=
\int_{0}^{+\infty}\int_{\mathbb{R}^{3}}
2\mu e^{it\mu^{2}}
\Im
\left(
e^{-i\mu |x-y_{j}|}
e^{i \mu |y-y_{\ell}|}
c_{j\ell}(\mu)
\right)
\frac{f(y)}{|x-y_{j}||y-y_{\ell}|}
\psi(\mu/M)
dyd\mu
\end{equation}
If we can prove a dispersive estimate
for $I_{M}$ which is uniform in $M$, this will also give an estimate
for $I$. Since $c_{j\ell}(\mu)$ is bounded near $\mu=0$
by assumption and $c_{j\ell}(0)$ is real, while
\begin{equation*}
\Im( e^{-i\mu |x-y_{j}|}e^{i\mu |y-y_{\ell}|})
=\sin(\mu ( |x-y_{\ell}|-|y-y_{j}|))
\end{equation*}
vanishes for $\mu=0$, we can integrate by parts with respect to
$\mu$ and we obtain
\begin{equation}\label{eq.IMter}
I_{M}=
\frac i t
\int_{0}^{+\infty}\int_{\mathbb{R}^{3}}
e^{it\mu^{2}}
\partial_{\mu}
\Im
\left(
e^{-i\mu |x-y_{j}|}
e^{i \mu |y-y_{\ell}|}
c_{j\ell}(\mu)
\psi(\mu/M)
\right)
\frac{f(y)}{|x-y_{j}||y-y_{\ell}|}
dyd\mu.
\end{equation}
For reasonable data $f$ the two integrals can be swapped using
Fubini's theorem and we arrive at
\begin{equation}\label{eq.IMqu}
I_{M}=
\frac i t
\int_{\mathbb{R}^{3}}
\frac{f(y)}{|x-y_{j}||y-y_{\ell}|}
\int_{0}^{+\infty}
e^{it\mu^{2}}
\partial_{\mu}
\Im
\left(
e^{-i\mu |x-y_{j}|}
e^{i \mu |y-y_{\ell}|}
c_{j\ell}(\mu)
\psi(\mu/M)
\right)
d\mu dy.
\end{equation}
The core of our proof will be a suitable estimate of the inner integral in the
variable $\mu$; expanding the derivative we obtain three
terms
\begin{equation*}
\int_{0}^{+\infty}
e^{it\mu^{2}}
\partial_{\mu}
\Im
\left(
e^{-i\mu |x-y_{j}|}
e^{i \mu |y-y_{\ell}|}
c_{j\ell}(\mu)
\psi(\mu/M)
\right)
d\mu=I_{1}+I_{2}+I_{3}
\end{equation*}
and precisely
\begin{equation}\label{eq.I1}
I_{1}=
(|y-y_{\ell}|- |x-y_{j}|)
\int_{0}^{+\infty}
e^{it\mu^{2}}
\Re
\left(
e^{i \mu |y-y_{\ell}|-i\mu |x-y_{j}|}
c_{j\ell}(\mu)
\psi(\mu/M)
\right)
d\mu,
\end{equation}
\begin{equation}\label{eq.I2}
I_{2}=
\int_{0}^{+\infty}
e^{it\mu^{2}}
\Im
\left(
e^{-i\mu |x-y_{j}|}
e^{i \mu |y-y_{\ell}|}
c_{j\ell}'(\mu)
\psi(\mu/M)
\right)
d\mu,
\end{equation}
and
\begin{equation}\label{eq.I3}
I_{3}=
\int_{0}^{+\infty}
e^{it\mu^{2}}
\Im
\left(
e^{-i\mu |x-y_{j}|}
e^{i \mu |y-y_{\ell}|}
c_{j\ell}(\mu)
\psi'(\mu/M)
\frac1M
\right)
d\mu.
\end{equation}

Consider the first integral $I_{1}$, or rather the integral
\begin{equation}\label{eq.I11}
\widetilde I_{1}=
\int_{0}^{+\infty}
e^{it\mu^{2}}
\Re
\left(
e^{i \mu |y-y_{\ell}|-i\mu |x-y_{j}|}
c_{j\ell}(\mu)
\psi(\mu/M)
\right)
d\mu.
\end{equation}
Inspired by a clever idea of Rodnianski and Schlag (\cite{RodnianskiSchlag04-disp}), we remark that
an integral of the form
\begin{equation*}
u(t,A)=
\int_{0}^{+\infty}
e^{it\mu^{2}}
\Re
\left(
e^{i A\mu}
c_{j\ell}(\mu)
\psi(\mu/M)
\right)
d\mu
\end{equation*}
is a solution of the one dimensional Schr\"odinger equation
$i\partial_{t}u=\partial^{2}_{A}u$. Thus we can apply the classical dispersive
estimate
\begin{equation*}
|u(t,A)|\leq C t^{-1/2}\int|u(0,A)|dA
\end{equation*}
which in our case gives
\begin{equation}\label{eq.part1}
\left|
\widetilde I_{1}
\right|
\leq \frac C{\sqrt t}\cdot \int|J_{1}(A)|dA
\end{equation}
with
\begin{equation}\label{eq.J1}
J_{1}(A)=
\int_{0}^{+\infty}
\Re
\left(
e^{i A\mu}
c_{j\ell}(\mu)
\psi(\mu/M)
\right)
d\mu.
\end{equation}
The same argument, used for $I_{3}$, shows that
\begin{equation}\label{eq.part3}
\left|
I_{3}
\right|
\leq \frac C{\sqrt t}\cdot \int|J_{3}(A)|dA
\end{equation}
where
\begin{equation}\label{eq.J3}
J_{3}(A)=
\int_{0}^{+\infty}
\Im
\left(
e^{i \mu A}
c_{j\ell}(\mu)
\psi'(\mu/M)
\frac1M
\right)
d\mu.
\end{equation}
A similar argument can be applied to $I_{2}$. First of all we
recall that the matrix $\Gamma_{\alpha,Y}(\mu)$
for $\mu\in[0,+\infty)$ has the form
(see \eqref{eq.Gamma1})
\begin{equation}\label{eq.gam3}
\Gamma_{\alpha,Y}(\mu)=-\frac{i\mu}{4\pi}I
+\mathrm{diag}[\alpha_{1},\dots,\alpha_{N}]+A(\mu)
\end{equation}
where $A(\mu)$ has coefficients
\begin{equation}\label{eq.Amu}
[A(\mu)]_{jj}=0,\qquad
[A(\mu)]_{j\ell}=-
\frac{e^{i\mu|y_{j}-y_{\ell}|}}{4\pi|y_{j}-y_{\ell}|}\quad
\text{for $j\neq\ell$.}
\end{equation}
Hence the derivative $\partial_{\mu}\Gamma_{\alpha,Y}(\mu)$
is simply
\begin{equation}\label{eq.derG}
[\partial_{\mu}\Gamma_{\alpha,Y}(\mu)]_{j\ell}=
-\frac i{4\pi} e^{i\mu|y_{j}-y_{\ell}|}.
\end{equation}
Since
\begin{equation*}
\partial_{\mu}[\Gamma_{\alpha,Y}(\mu)^{-1}]=
\Gamma_{\alpha,Y}(\mu)^{-1}
(\partial_{\mu}\Gamma_{\alpha,Y}(\mu))
\Gamma_{\alpha,Y}(\mu)^{-1},\qquad
c_{j\ell}(\mu)=
[\Gamma_{\alpha,Y}(\mu)^{-1}]_{j\ell},
\end{equation*}
we obtain the representation
\begin{equation}\label{eq.derc}
c'_{j\ell}(\mu)=
-\frac i{4\pi}
\sum_{k_{1},k_{2}=1}^{N}
c_{jk_{1}}(\mu)
e^{i\mu|y_{k_{1}}-y_{k_{2}}|}
c_{k_{2}\ell}(\mu)
\end{equation}
If we plug this into $I_{2}$ we obtain
\begin{equation}\label{eq.I22}
I_{2}= \frac1{4\pi}\sum_{k_{1},k_{2}=1}^{N}
\int_{0}^{+\infty}
e^{it\mu^{2}}
\Re
\left(
c_{jk_{1}}(\mu)c_{k_{2}\ell}(\mu)
e^{i \mu |y-y_{\ell}|-i\mu |x-y_{j}|
+i\mu|y_{k_{1}}-y_{k_{2}}|}
\psi(\mu/M)
\right)
d\mu.
\end{equation}
The same argument used above for $I_{1}$ and $I_{3}$ gives that
\begin{equation}\label{eq.part2}
\left|
I_{2}
\right|
\leq \frac C{\sqrt t}\cdot \int|J_{2}(A)|dA
\end{equation}
where
\begin{equation}\label{eq.J2}
J_{2}(A)=
\int_{0}^{+\infty}
\Re
\left(
e^{i \mu A}
c_{jk_{1}}(\mu)c_{k_{2}\ell}(\mu)
\psi(\mu/M)
\right)
d\mu.
\end{equation}

We shall now prove
the estimate
\begin{equation}\label{eq.goal}
\int(|J_{1}(A)|+|J_{2}(A)|+|J_{3}(A)|)dA<C
\end{equation}
for some constant independent of $M$, from which the Theorem
will follow immediately. To this end we shall need to study
the coefficients $c_{j\ell}(\mu)$ closer.

\begin{lemma}\label{lem.cjl}
Let $\Gamma_{\alpha,Y}(\mu)$ be the matrix defined in \eqref{eq.gam3},
\eqref{eq.Amu} and assume it is invertible with inverse locally bounded
for $\mu\in[0,+\infty)$.

Then the coefficients $c_{j\ell}(\mu)$
of the inverse matrix $\Gamma_{\alpha,Y}(\mu)^{-1}$
satisfy the following properties:

(i) The coefficients are holomorphic on a neighbourhood of
the positive real axis; their derivatives satisfy the estimates
\begin{equation}\label{eq.estc}
|c_{j\ell}(\mu)|\leq C\mm^{-1},\qquad
|c'_{j\ell}(\mu)|\leq C\mm^{-2},\qquad
|c''_{j\ell}(\mu)|\leq C\mm^{-2},\qquad
\mu\geq0
\end{equation}
where we used the notation $\mm=(1+\mu^{2})^{1/2}$.

(ii) The coefficients can be written as
\begin{equation}\label{eq.asy}
c_{j\ell}(\mu)=
{4\pi i\delta_{j\ell}}
{\mm}^{-1}+d_{j\ell}(\mu),
\end{equation}
or equivalently
\begin{equation}\label{eq.asy2}
c_{j\ell}(\mu)=
{4\pi i\delta_{j\ell}}
\mu{\mm}^{-2}+\tilde d_{j\ell}(\mu),
\end{equation}
where the holomorphic functions $d_{j\ell},\tilde d_{j,\ell}$
satisfy the inequalities
\begin{equation}\label{eq.djl}
|d_{j\ell}|+|d'_{j\ell}|+|d''_{j\ell}|
\leq C\mm^{-2},\qquad
| \tilde d_{j\ell}|+| \tilde d'_{j\ell}|+| \tilde d''_{j\ell}|
\leq C\mm^{-2},\qquad\mu\geq0.
\end{equation}
\end{lemma}

\begin{proof}
It is clear from the definition that $c_{j\ell}(\mu)$ are holomorphic
(as soon as they are defined); moreover, by formula \eqref{eq.derc}
proved above we see that the second and third estimates
in \eqref{eq.estc} are immediate
consequences of the first one in \eqref{eq.estc}.

In order to prove the first inequality in \eqref{eq.estc},
it is sufficient to recall formulas \eqref{eq.gam3},
\eqref{eq.Amu}, i.e.,
\begin{equation*}
\Gamma_{\alpha,Y}(\mu)=C_{0}\mu I+D+A(\mu),
\end{equation*}
with $C_{0}=-i/4\pi$, $D$ a constant diagonal matrix, and $A(\mu)$
a matrix with bounded coefficients. From this expression it is clear
that the entries of the inverse matrix $\Gamma_{\alpha,Y}(\mu)^{-1}$
are of the form
\begin{equation*}
c_{j\ell}(\mu)=
\frac{\pm (C_{0}\mu)^{N-1}+P(\mu)}
{\pm (C_{0}\mu)^{N}+Q(\mu)}
\end{equation*}
where $P$, $Q$ are functions of $\mu$ of order $N-2$ and $N-1$ respectively:
\begin{equation*}
    |P(\mu)|\leq C\mm^{N-2},\qquad
    |Q(\mu)|\leq C\mm^{N-1}.
\end{equation*}
Taking into account the
assumption that $c_{j\ell}(\mu)$ are locally bounded functions,
estimate \eqref{eq.estc} follows easily.

In order to prove the asymptotic
expansion \eqref{eq.asy}, we write
\begin{equation*}
I=\Gamma_{\alpha,Y}(\mu)^{-1}\Gamma_{\alpha,Y}(\mu)=
\Gamma_{\alpha,Y}(\mu)^{-1}(C_{0}\mu I+D+A(\mu))=
C_{0}\mu\Gamma_{\alpha,Y}(\mu)^{-1}+ \Gamma_{\alpha,Y}(\mu)^{-1}[D+A(\mu)]
\end{equation*}
which implies ($C_{0}^{-1}=4\pi i$)
\begin{equation*}
\Gamma_{\alpha,Y}(\mu)^{-1}=
4\pi i\mu^{-1}
I-
[D+A(\mu)]\cdot
\Gamma_{\alpha,Y}(\mu)^{-1}
4\pi i\mu^{-1}
\end{equation*}
(it is sufficient to prove the estimates for $\mu>1$ since we already know
that the functions are smooth near 0).
Notice that the last term is bounded by $C\mm^{-2}$.
Now, using the elementary identity
\begin{equation*}
\frac1\mu=\frac1{(1+\mu^{2})^{1/2}}+
\frac1{\mu(1+\mu^{2})^{1/2}[\mu+(1+\mu^{2})^{1/2}]}
\end{equation*}
and estimate \eqref{eq.estc} already proved, we
easily obtain \eqref{eq.asy} and the estimate
\begin{equation*}
|d_{j\ell}(\mu)|\leq C\mm^{-2}.
\end{equation*}
The remaining estimates on the derivatives of $d_{j\ell}$
follow immediately from \eqref{eq.asy} and \eqref{eq.estc}.

The second expansion \eqref{eq.asy2} is proved in an identical way; indeed,
\begin{equation*}
\frac1\mu=\frac\mu{1+\mu^{2}}+
\frac1{\mu(1+\mu^{2})}.
\end{equation*}
\end{proof}

\begin{remark}\label{rem.odd}
The reason for the two different asymptotic expansions in the Lemma
is the following: the main term decays like $\mu^{-1}$ and must
be treated with care, using explicit formulas for the Fourier transform.
Since the Fourier transform here appears as a sinus or cosinus
transform only, it is quite useful to have its expression
both as an \emph{even} and
an \emph{odd} function of $\mu$. The remainder terms are harmless
since they give rise to integrals which can be easily estimated uniformly in $M$.
\end{remark}

We are ready to estimate the $L^{1}$ norms on $\mathbb{R}$
of the three quantities $J_{1}(A),J_{2}(A),J_{3}(A)$.
Consider the quantity $J_{1}(A)$; using the expansion \eqref{eq.asy2}
we have
\begin{equation}\label{eq.split}
J_{1}(A)=J^{a}_{1}(A)+J_{1}^{b}(A)
\end{equation}
where
\begin{equation*}
J^{a}_{1}(A)=
4\pi\int_{0}^{+\infty}
\Re
\left(
e^{i A\mu}4\pi i
\mu\mm^{-2}
\psi(\mu/M)
\right)
d\mu=
-(4\pi)^{2}
\int_{0}^{+\infty}
\sin({A\mu})
\mu\mm^{-2}
\psi(\mu/M)
d\mu
\end{equation*}
and
\begin{equation*}
J^{b}_{1}(A)=
\int_{0}^{+\infty}
\Re
\left(
e^{i A\mu}
\tilde d_{j\ell}(\mu)
\psi(\mu/M)
\right)
d\mu.
\end{equation*}
The remainder term $J^{b}(A)$ can be estimated directly.
Indeed, by \eqref{eq.djl} it is easy to see that $ J^{b}_{1}(A)$ is
bounded
\begin{equation*}
| J^{b}_{1}(A)|\leq
C\int_{0}^{+\infty}
\mm^{-2}d\mu\leq C
\end{equation*}
with $C$ independent of $M,A$.
Moreover, writing $e^{iA\mu}=(iA)^{-1}\partial_{\mu}e^{iA\mu}$
we have
\begin{equation*}
J_{1}^{b}(A)=\frac{1} A
\int_{0}^{+\infty}
\Im
\left(
\partial_{\mu}
e^{i A\mu}
\tilde d_{j\ell}(\mu)
\psi(\mu/M)
\right)
d\mu
\end{equation*}
but $c_{j\ell}(0)$ is real, so integrating by parts we get
\begin{equation*}
J_{1}^{b}(A)=-\frac1 A
\int_{0}^{+\infty}
\Im
\left(
e^{i A\mu}
\tilde d'_{j\ell}(\mu)
\psi(\mu/M)
\right)
d\mu
-\frac1 A
\int_{0}^{+\infty}
\Im
\left(
e^{i A\mu}
\tilde d_{j\ell}(\mu)
\psi'(\mu/M)
\frac1M
\right)
d\mu.
\end{equation*}
As before, using \eqref{eq.djl} we see that
both integrals are bounded by a constant
$C'$ independent of $M,A$ and this implies
\begin{equation*}
|AJ^{b}_{1}(A)|\leq C'.
\end{equation*}
Finally, we write $e^{iA\mu}=(iA)^{-1}\partial_{\mu}e^{iA\mu}$
one more time in the last formula:
\begin{align*}
J_{1}^{b}(A)=\frac{1} {A^{2}}
\int_{0}^{+\infty}
\Re
\left(
\partial_{\mu}
e^{i A\mu}
\tilde d'_{j\ell}(\mu)
\psi(\mu/M)
\right)
d\mu
+
\frac1 {A^{2}}
\int_{0}^{+\infty}
\Re
\left(
\partial_{\mu}
e^{i A\mu}
\tilde d_{j\ell}(\mu)
\psi'(\mu/M)
\frac1M
\right)
d\mu.
\end{align*}
We integrate again by parts; this time we obtain a boundary term
\begin{equation}\label{eq.bdr}
\frac1{A^{2}}
\Re \tilde d'_{j\ell}(0)\psi(0)+
\frac1{A^{2}}
\Re \tilde d_{j\ell}(0)\psi'(0)\frac1M
\end{equation}
and a few integrals
\begin{align}\label{eq.intgr}
\frac1 {A^{2}}
\int_{0}^{+\infty}
\Re
\left(
e^{i A\mu}
\tilde d''_{j\ell}(\mu)
\psi(\mu/M)
\right)
d\mu
+&
\frac2 {A^{2}}
\int_{0}^{+\infty}
\Re
\left(
e^{i A\mu}
\tilde d'_{j\ell}(\mu)
\psi'(\mu/M)
\frac1M
\right)
d\mu\\
+&
\frac1 {A^{2}}
\int_{0}^{+\infty}
\Re
\left(
e^{i A\mu}
\tilde d_{j\ell}(\mu)
\psi''(\mu/M)
\frac1{M^{2}}
\right)
d\mu.
\end{align}
Applying as above the estimates \eqref{eq.djl} of the Lemma we obtain
\begin{equation*}
|A^{2}J^{b}_{1}(A)|\leq C''
\end{equation*}
for a constant independent of $M,A$. Then we conclude
\begin{equation*}
\int|J_{1}^{b}(A)|dA\leq \int C(1+A^{2})^{-1}dA\leq C
\end{equation*}
for some bound $C$ uniform in $M$.

Let us now consider
the main term $J^{a}_{1}(A)$ in \eqref{eq.split}. It
is a standard sinus-transform, which can be rewritten
as a Fourier transform as follows (recall $\psi$ is an even function)
\begin{equation*}
J^{a}_{1}(A)=
8\pi^{2}i
\int_{-\infty}^{+\infty}
e^{iA\mu}
\frac\mu{1+\mu^{2}}\psi(\mu/M)d\mu=
C\cdot
\left( \mathcal{F}\left(g_{1}\right)
\star
\mathcal{F}\left(g_{2}\right) \right)(A)
\end{equation*}
where we have used the convolution theorem and we have denoted

\begin{equation}
g_{1}(\mu)=\frac{\mu}{1+ \mu^{2}}, \;\;\;\;\; g_{2}(\mu)=
\psi(\mu / M)
\end{equation}

Since
\begin{equation*}
\left(\mathcal{F}(g_{1})\right)(A)
=C \frac A{|A|}e^{-|A|}
\end{equation*}
and
\begin{equation*}
\left(\mathcal{F}(g_{2})\right)(A)=
M(\mathcal{F}\psi) (AM)
\end{equation*}
we obtain
\begin{equation*}
\|J^{a}_{1}\|_{L^{1}(\mathbb{R})}\leq
C\left\|\mathcal{F}(g_{1})
\right\|_{L^{1}(\mathbb{R})}
\|\mathcal{F}(g_{2})\|_{L^{1}(\mathbb{R})} \leq C
\end{equation*}
for some constant independent of $M$. In conclusion we have proved the
bound $\|J_{1}\|_{L^{1}(\mathbb{R})}<C$, uniformly in $M$.

The corresponding estimate for $J_{3}(A)$, see \eqref{eq.J3}, is
almost identical; the only modification is
the use of the explicit formula for the Fourier transform of $\mm^{-1}$,
which is
\begin{equation*}
\mathcal{F}(\mm^{-1})
=C K_{0}\left(|A|\right).
\end{equation*}
Here $K_{0}(s)$ is the modified Bessel function of order 0, whose behaviour is
the following:
\begin{equation*}
|K_{0}(s)|\leq C|\log s|,\quad 0<s\leq 1;\qquad
|K_{0}(s)|\leq Cs^{-\frac12}e^{-s},\quad s\geq 1.
\end{equation*}
In particular, it is clear that $ \mathcal{F}(\mm^{-1})$ belongs
to $L^{1}(\mathbb{R})$ which is what is needed in our
computation.
On the other hand the estimate for $J_{2}(A)$ (see \eqref{eq.J2}) is
immediate since the integrand decays as $\mm^{-2}$.
This concludes the proof of \eqref{eq.goal}.

Recalling now
\eqref{eq.I1}, \eqref{eq.I2}, \eqref{eq.I3} and
\eqref{eq.part1}, \eqref{eq.part2},\eqref{eq.part3}
we have
\begin{equation*}
|I_{1}|\leq\frac C{t^{1/2}}(|y-y_{\ell}|+|x-y_{j}|),\qquad
|I_{2}|+|I_{3}|\leq \frac C{t^{1/2}}
\end{equation*}
with $C$ independent of $M$ and of $x$. Thus,
recalling \eqref{eq.IMter} and taking the
limit as $M\to\infty$, we obtain the following
estimate for the integral $I$ in \eqref{eq.I}:
\begin{equation*}
|I|\leq \frac C{t^{3/2}}
\int_{\mathbb{R}^{3}}
|f(y)|\frac{1+|y-y_{\ell}|+|x-y_{j}|}{|y-y_{\ell}||x-y_{j}|}
dy
\end{equation*}
Since $P_{ac}e^{itH_{\alpha,Y}}f$ is the sum of the integrals
$I$ for all couples $j,\ell=1,\dots,N$, plus the
free wave $e^{-it\Delta}f$, we obtain
\begin{equation*}
|(P_{ac}e^{itH_{\alpha,Y}}f)(x)|\leq
\frac C{t^{3/2}}
\int
|f(y)|
\sum_{j,\ell =1}^{N}
\left(
1+\frac{1+|y-y_{\ell}|+|x-y_{j}|}{|y-y_{\ell}||x-y_{j}|}
\right)dy.
\end{equation*}
The inner sum gives
\begin{equation*}
\sum_{j,\ell=1}^{N}
\left(
1+\frac{1+|y-y_\ell|+|x-y_j|}
{|y-y_\ell||x-y_j|}
\right)=
\sum_{j,\ell =1}^{N}
\left(
1+\frac1{|y-y_\ell|}
\right)
\cdot
\left(
1+\frac1{|x-y_j|}
\right)
=
w(x)w(y)
\end{equation*}
where $w(x)$ is the weight function \eqref{eq.weight}.

In conclusion we have proved
\begin{equation*}
w(x)^{-1}
|(P_{ac}e^{itH_{\alpha,Y}}f)(x)|\leq
\frac C{t^{3/2}}
\int
|f(y)| w(y)dy
\end{equation*}
as claimed.

\vspace{1cm}

\bibliographystyle{plain}

\vspace{1cm}
\end{document}